\documentclass[leqno,12pt]{article}
\usepackage{graphicx,fancybox,ascmac}
\usepackage{amsmath, amssymb}
\numberwithin{equation}{section}

\newtheorem{theorem}{Theorem}[section]
\newtheorem{cor}{Corollary}[section]
\newtheorem{prop}{Proposition}[section]
\newtheorem{ex}{Example}[section]
\newtheorem{lem}{Lemma}[section]
\newtheorem{rem}{Remark}[section]

\numberwithin{equation}{section}

\newcommand{\nn}{\nonumber}
\newcommand{\edn}{\end{description}}
\newcommand{\en}{\begin{enumerate}}
\newcommand{\een}{\end{enumerate}}
\begin{document}

\begin{center}
{\bf \Large
Binomial series and complex difference equations
}\vspace{0.5325cm}

{Katsuya Ishizaki\footnote{\ Supported by the discretionary budget (2018) of the President of the Open University of Japan.}
 and Zhi-Tao Wen\footnote{\ was supported by the National Natural Science Foundation of China (No.~11971288 and No.~11771090) and Shantou University SRFT (NTF18029).}
 }
\vspace{0.125cm}

{\it \small $^{*1}$ The Open University of Japan, 2-11 Wakaba, Mihama-ku, Chiba, 261- 8586 Japan, ishizaki@ouj.ac.jp}

{\it \small $^{*2}$ Department of Mathematics, Shantou University, Daxue Road No.~243, Shantou 515063, Guangdong, China, zhtwen@stu.edu}
\vspace{0.25cm}
\end{center}

\noindent {\bf Abstract.}\enspace We consider properties of binomial series
$\sum_{n=0}^\infty a_n z^{\underline{n}}$, where $z^{\underline{n}}=z(z-1)\cdots(z-n+1)$ and the convergence of binomial series in the complex domain.
The order of growth of entire and meromorphic solutions of some difference equations represented by binomial series are discussed. Examples are given. As an application, we construct a difference Riccati equation possessing a transcendental meromorphic solution of order $1/2$.
\vspace{0.125cm}

\begin{flushleft}
\end{flushleft}

\begin{flushleft}
\end{flushleft}

\section{Introduction}
We consider properties of binomial series in the complex domain and discuss a criteria for convergence of binomial series in connection with the order of growth of entire functions.
The motivations of our concern are from the study of integrability of discrete functional equations in the complex domain, which has been actively discussed by several detectors, see e.g.,~\cite{AHH2000},~\cite{GHRV2009},~\cite{H2017},~\cite{HK2006}. Once a candidate of integrable discrete equation is chosen by a detector to check the necessary condition of existence of solutions, it is required to construct a rigorous solution.
The binomial series has been expected for the construction of analytic functions in the complex domain. In this paper, we recall accumulated researches of binomial series and give an idea to construct entire and meromorphic solutions of some difference equations in the complex domain.

For a function $f$, we denote by $\Delta f(z)=f(z+1)-f(z)$ the difference operator. Let $n$ be a nonnegative integer. Define $\Delta^n f(z)=\Delta(\Delta^{n-1} f(z))$ for $n\geq 1$, and write $\Delta^0 f=f$.
Define $z^{\underline{0}}=1$ and
\begin{equation}
z^{\underline{n}}=z(z-1)\cdots(z-n+1)=n!\begin{pmatrix}
z\\
n
\end{pmatrix},\quad n=1, 2, 3, \dots,\label{1.3}
\end{equation}
which is called a {\it falling factorial}, see Subsection~\ref{subH},

Consider the formal series of the form
\begin{equation}
Y(z)=\sum_{n=0}^\infty a_n z^{\underline{n}},\quad a_n\in\mathbb C,\quad n=0,1,2, \dots.\label{1.5}
\end{equation}
For a fixed $z$, if $\sum_{n=0}^\infty |a_n||z^{\underline{n}}|$ converges, we say $Y(z)$ in \eqref{1.5} absolutely converges at $z$. We write  the limit function of $Y(z)$ as $y(z)$.
Let $\{\alpha_n\}$ be a sequence. We define a quantity concerning  $\{\alpha_n\}$ as
\begin{equation}
\chi(\{\alpha_n\})=\limsup_{n\to\infty}\frac{n\log n}{-\log|\alpha_n|}.\label{3.1}
\end{equation}

One of the main results in this paper is the following theorem.
\begin{theorem} \label{Convergence}
Suppose that $\chi(\{a_n\})<1$. Then the formal series $Y(z)$ given by \eqref{1.5} converges to $y(z)$ uniformly on every compact subset in $\mathbb C$. Moreover, the order of growth of $y(z)$ coincides with $\chi(\{a_n\})$.
\end{theorem}
We recall in Section~\ref{Backgrounds} backgrounds of researches of difference calculus and difference equations in the complex domain, and also recall the theory of entire functions related to the order of growth and interpolations.
In Section~\ref{Preliminaries}, we show some preliminary lemmas for the proof of Theorem~\ref{Convergence}.
We give the proof of Theorem~\ref{Convergence} in Section~\ref{SConvergence} and add further results on the case when $\chi(\{a_n\})=1$ with some conditions.
In Section~\ref{examples}, we give examples of formal solutions of some linear difference equations by binomial series and discuss the convergence of them.
Section~\ref{Growth} devotes to the observations of the order of growth of linear difference equations with polynomial coefficients and discussions of a difference Riccati equation with a rational function as coefficients by means of Nevanlinna theory.

\section{Backgrounds of difference calculus and entire functions}\label{Backgrounds}
In this section, we recall backgrounds of researches of difference equations intersecting the complex analysis. By~\cite{AHH2000}, Yanagihara's works in 1980's concerning complex difference equations, e.g., \cite{Yanagihara1980},~\cite{Yanagihara1985}, were reviewed.
Then the results in early 2000's on difference analogues of the value distribution theory, e.g.,~\cite{CF2008},~\cite{HKT2014}, generalized Yanagihara's theorems with some conditions on the order of growth.

We further look back to contacts of binomial series in difference calculus and the value distribution theory on entire functions.  In Subsection~\ref{subH}, we recall the falling factorial in difference calculus, and Striling's numbers in Subsection~\ref{subSt}, which had been studied before the middle of 19 century when complex analysis  was actively investigated.
In Subsections~\ref{subBS} and ~\ref{subCB}, we recall the researches on binomial series and its applications to linear difference equations due to N\"{o}rlund, Milne-Thomson and others.
We again go back to the late of of 19th century and the early of 1900's in Subsection~\ref{subOG} to see the
the order of growth of entire functions. These works lead to the Carlson theorem in Subsection~\ref{subCN} concerning an analytic function and binomial series.

\subsection{Difference calculus}\label{subH}
We back to Newton's Principia~\cite{Newton1687} in 1687. Indeed, we will discuss Newton series below which probably named after his name. Boole~\cite[Page 61]{Boole1860} mentioned that Newton's Principia is the first attempt at finding a general interpolation formula.
According to Jordan~\cite{Jordan1950}, Taylor ascribed the difference calculus~\cite{Taylor1715} in 1715 and Stirling established the theory of the difference calculus and gave useful methods introducing the  Stirling numbers in~\cite{Stirling1730} in 1730.
It is also mentioned in~\cite{Jordan1950} that Euler in 1755 introduced the symbol $\Delta$ for the differences in~\cite{Euler1755}, which is universally used now.
We here mention that the list of references in N\"{o}rlund's book~\cite[Pages 244--531]{Norlund1924} is very rich and is very useful when looking up classic literature.

We defined the falling factorial by \eqref{1.3}, which yields
\begin{equation}
\Delta z^{\underline{n}}=(z+1)^{\underline{n}}-z^{\underline{n}}=nz^{\underline{n-1}}\ .\label{21.02}
\end{equation}
This corresponds to $(z^n)'=nz^{n-1}$ in the differential calculus.
The symbols of the falling factorial are variously different for each researcher, but they gave the same comment that \eqref{21.02} holds.
For examples, Boole~\cite[Page 6]{Boole1860} and Milne-Thomson~\cite[Page 25]{M-Thomson1933} used the symbol $z^{(n)}$ to represent \eqref{1.3} and called it factorial expression.
Jordan adopted the symbol $(z)_n$ and called it factorial of degree $n$ in~\cite[Page 45]{Jordan1950}.
N\"{o}rlund used $z(z-1)\cdots(z-n+1)$ as it is and mentioned \eqref{21.02}~in \cite[Page 5, (9)]{Norlund1924}. In this paper, we adopt \eqref{1.3}, following \cite[Page 10]{Aigner2007} and \cite[Page 17]{KP2001}.
The formula \eqref{21.02} is one of the key idea when we consider the connection between difference calculus and differential calculus.

We call a function $f(z+n)$ an {\it $n$-th  shift} of $f(z)$ for $n\in\mathbb N$, in particular, we call $f(z+1)$ a {\it shift} of $f(z)$ simply for $n=1$.
We have relations between higher order differences and higher order shifts
\begin{equation}
\Delta^nf(z)=\sum_{j=0}^n \begin{pmatrix} n  \\ j \end{pmatrix}(-1)^{n-j}f(z+j),\label{1.6}
\end{equation}
and
\begin{equation}
f(z+n)=\sum_{j=0}^n \begin{pmatrix} n  \\ j \end{pmatrix}\Delta^jf(z).
\label{1.7}
\end{equation}
The relations \eqref{1.6} and \eqref{1.7} are stated in e.g.,~\cite[Page 4, (5), (7)]{Norlund1924}.
By means of \eqref{1.6} and \eqref{1.7}, it is possible that we represent a difference equation by $\Delta$ and also by shifts.

\subsection{The Stirling numbers}\label{subSt}

We next consider the Stirling numbers, which connect the falling factorial given by \eqref{1.3} and the monomial $z^n$.
Let $n\geq2$ be an integer. Write
\begin{equation}
z^{\underline{n}}=z^n+\sum_{j=1}^{n-1}\eta_{j, n}\, z^j,\label{21.03}
\end{equation}
where $\eta_{n,n}=1$ and $\eta_{0,n}=0$. It is known that $\eta_{j,n}$ satisfies a recurrence relation
\begin{equation}
\eta_{j, n+1}=\eta_{j-1,n}-n\, \eta_{j,n},\quad 1\leq j\leq n.\label{21.04}
\end{equation}
We also write
  \begin{equation}
z^{n}=z^{\underline{n}}+\sum_{k=1}^{n-1}\widetilde{\eta}_{k, n}\, z^{\underline{k}},\label{21.05}
    \end{equation}
where $\widetilde{\eta}_{n,n}=1$ and $\widetilde{\eta}_{0,n}=0$.
Then we have
    \begin{equation}
\widetilde{\eta}_{k, n+1}=\widetilde{\eta}_{k-1,n}+k\,\widetilde{\eta}_{k,n},\quad 1\leq k\leq n.\label{21.06}
    \end{equation}
The coefficients $\eta_{j,n}$ in \eqref{21.03} and $\tilde\eta_{k,n}$ in \eqref{21.05} are called the first kind and the second kind of Stirling numbers for the falling factorial   $z^{\underline{n}}$, respectively. The recurrence relations \eqref{21.04} and \eqref{21.06} are known as the Pascal rule for Stirling numbers, see e.g.,~\cite[Chapter 1]{Aigner2007}.
The symbols of the Stirling numbers are variously different for each researcher
as we saw the case of the falling factorial.
For examples, N\"{o}rlund indeed discussed relations between polynomials and binomial (factorial) polynomials in \cite[Page 148]{Norlund1924} using the symbol
$(-1)^j\binom{n}{j}B_j^n$ corresponding to the first kind Stirling number.
Milne-Thomson~\cite[Chapter~6]{M-Thomson1933} studied an expansion of $z^{\underline{n}}$ in powers of $z$ with the same motivation.
Stirling numbers are considered deeply with corresponding relations to \eqref{21.03}--\eqref{21.06} in~\cite[Pages 142--224]{Jordan1950}, in which Jordan used $S_n^j$ for the first kind and $\mathfrak{S}_n^j$ for the second kind.

\subsection{Binomial series (Newton series, factorial series)}\label{subBS}
We are concerned with the formal series of the form \eqref{1.5} in this subsection and in the next subsection.
We call the right hand side of \eqref{1.5} the {\it binomial series} in this paper, see e.g.,~\cite[Page~198]{Hille1959},~\cite{IY2004}.
Historically, \eqref{1.5} has been variously called the Newton series or the factorial series.
For examples, Boole~\cite[Page 11]{Boole1860} wrote
$\phi(z)=a+bz+cz^{(2)}+dz^{(3)}+\cdots+hz^{(m)}+\cdots$ and called it a series of factorials, where Boole's notation $z^{(m)}$ is the falling factorial.
In~\cite[Page 15, (35)]{Norlund1924}, N\"{o}rlund treated series of the form
\begin{equation}
F_a(z)=\sum_{n=0}^\infty \frac{\Delta^n f(a)}{n!}(z-a)^{\underline{n}}\ ,\label{22.02}
\end{equation}
and Jordan also discussed \eqref{21.02} in~\cite[Pages 357--368]{Jordan1950}.
Further, N\"{o}rlund considerd~\cite[Pages 205, 222--233, (17)]{Norlund1924},
\begin{equation}
\tilde F(z)=\sum_{n=0}^\infty \tilde{a}_n \begin{pmatrix}
z-1\\ n
\end{pmatrix}=\sum_{n=0}^\infty \frac{\tilde{a}_n}{n!}(z-1)^{\underline{n}},
\quad \tilde{a}_n\in\mathbb C,\quad n=0,1, 2, \dots,\label{22.03}
\end{equation}
which was called Newton series by N\"{o}rlund as a general form of \eqref{22.02} when $a=1$.
Milne-Thomson~\cite[Page 57--60]{M-Thomson1933} studied \eqref{22.02} as Newton's general interpolation formula.

\subsection{Convergence of binomial series}\label{subCB}
We recall classical results of convergence of binomial series in the complex domain.
Landau~\cite{Landau1906} investigated the region of convergence of binomial series and obtained fundamental results.
N\"{o}rlund~\cite[Pages 205, 257--262]{Norlund1924} and Milne-Thomson~\cite[Chapter~10]{M-Thomson1933} deeply studied convergence of factorial series of the form \eqref{22.03}.
Following~\cite{Norlund1924} and~\cite{M-Thomson1933}, we turn around the discussions of convergence of binomial series.
We here use a symbol $\Re z$ to represent the real part of $z$.
\begin{theorem}\label{M-Thomson1} {\rm{\bf(Landau, N\"{o}rlund, Milne-Thomson)}}\enspace
If \eqref{1.5} converges for $z=z_0$, then \eqref{1.5} converges for $z$ such that $\Re z>\Re z_0$.
\end{theorem}
Theorem~\ref{M-Thomson1} says that the region of convergence of binomial series is the right half plane. We define a real number $\lambda$ such that \eqref{1.5} converges for
$\Re z>\lambda+\varepsilon$, where $\varepsilon$ is positive and arbitrarily small, and diverges for $\Re z<\lambda-\varepsilon$, and call $\lambda$ the  {\it abscissa} of convergence of \eqref{1.5}.
If $\lambda=\infty$, the series \eqref{1.5} is everywhere divergent. If $\lambda=-\infty$, and the series \eqref{1.5} converges in the whole complex plane. Landau~\cite{Landau1906} found a criteria for abscissa of convergence $\lambda$ using a sum of coefficients of binomial series, see e.g., Milne-Thomson~\cite[Page~279]{M-Thomson1933}. We will see $\lambda=-\infty$ when $\chi(\{a_n\})<1$ and $\lambda\leq 0$ when $\chi(\{a_n\})=1$ with some conditions in Section~\ref{SConvergence}.

\subsection{The order of growth of entire functions}\label{subOG}

Liouville~\cite{Liouville1879} used a symbol $M(x,y)$, $z=x+y\sqrt{-1}$ to express the ``modulus'' of a function $f(z)$ in a certain $xy$ geometric rectangle to state the well known Liouville's theorem, which was obtained in 1847. In this connection, Lindel\"{o}f~\cite[Page 371]{Lindelof1905}, Phragm\'{e}n and Lindel\"{o}f~\cite{PL1908} and Wiman~\cite{Wiman1914} used a symbol $M(r)$ of maximum for $|z|=r$ for instance.
Constructing an entire function $f(z)$ with given value distribution was considered by Weierstrass~\cite{Weierstrass1876} in 1876. The Weierstrass factorization theorem asserts that every entire function can be represented as a product involving its zeroes, which is possibly infinite, see e.g.,~\cite[Pages 18--24]{Boas1954}, \cite[Page 225--229]{Hille1959}, \cite[Pages 282--285]{Markushevich1965}.
These are classical ideas in study of the growth and the value distribution of entire functions.

We now start with \cite[Pages 8--9]{Boas1954} to recall again the study of the growth of entire function $f(z)$. The {\it order of growth} $\rho=\rho(f)$ is defined by
\begin{equation}
\rho=\limsup_{r\to\infty}\frac{\log\log M(r,f)}{\log r},\quad 0\leq \rho\leq \infty.\label{32.01}
\end{equation}
When $0<\rho<\infty$, the type $\tau=\tau(f)$ is defined by
\begin{equation}
\tau=\limsup_{r\to\infty}\frac{\log M(r,f)}{r^{\rho}},\quad 0\leq \tau\leq \infty.
\label{32.02}
\end{equation}
It is known that the order $\rho$ is finite if and only if for every positive integer $\varepsilon$ but for no negative $\varepsilon$,
\begin{equation}
M(r,f)=O(e^{r^\rho+\varepsilon}),\quad\text{as $r\to\infty$}.\label{32.03}
\end{equation}
According to Nevanlinna~\cite[Page30]{Nevanlinna1929}, the order of growth was introduced by Borel~\cite[Page 362]{Borel18970} in 1897 adopting the method like \eqref{32.03}.
The concept of iterating order was also introduced in this Borel's paper. Pringsheim~\cite[Page 260--263]{Pringsheim1904} and Lindel\"{o}f~\cite[Pages 373--374]{Lindelof1905} also employed \eqref{32.03}.

We write $f(z)$ in the Taylor series
\begin{equation}
f(z)=\sum_{n=0}^\infty a_nz^n\label{32.04}
\end{equation}
and define
\begin{equation}
\chi=\limsup_{n\to\infty}\frac{n\log n}{-\log|a_n|}.\label{32.05}
\end{equation}
The relation between the order of growth of $f(z)$ and the coefficients of Taylor series \eqref{32.04} was discussed by Lindel\"{o}f~\cite[Chapter III]{Lindelof1902},~\cite{Lindelof1903} and deepened by Pringsheim~\cite[Page 260--263]{Pringsheim1904} comparing $\sqrt[n]{|a_n|}$ and $n^{-\rho\pm\varepsilon}$, where $\varepsilon$ is a small number. It is well known that  $\rho=\chi$, see e.g., \cite[Pages 9--11]{Boas1954},~\cite[Pages 24--25]{JV1985},~\cite[Pages 253--254]{Titchmarsh1939}. Following~\cite[Page 24]{JV1985}, we call \eqref{32.05} the Lindel\"of--Pringsheim formula, and we call $\rho=\chi$ the Lindel\"of--Pringsheim theorem in this paper.

\subsection{Carlson's results and Newton series}\label{subCN}

We recall a relation between Newton series and complex analysis by referring the classical results due to Carlson~\cite{Carlson1914}, see also e.g.,~\cite[Page 171]{Boas1954},~\cite{GR1978},~\cite{Hardy1920},~\cite[Page 58]{Leivn1996},~\cite{Pila2005},~\cite{Rubel1956}.
Main tools of the proof of the Carlson theorem are results obtained by Phragm\'{e}n--Lindel\"{o}f~\cite{PL1908}, see also e.g.,~\cite[Pages 37--39]{Leivn1996}.

Following~\cite[Pages 328--330]{Hardy1920},~\cite[Pages 185--186]{Titchmarsh1939}, we mention below the Carlson theorem and its Corollary as an application to Newton series \eqref{22.02}.

\begin{theorem}\label{ThmCarlson}{\rm\bf (Carlson)}
If \ {\rm(i)} $f(z)$ is regular at all points inside the angle $-\alpha<\theta<\alpha$, where $\alpha\geq \pi/2$, \ {\rm(ii)} \ $|f(z)|\leq Ae^{kr}$, where $k<\pi$, throughout this angle, {\rm(iii)} $f(n)=0$ for $n=1, 2, 3, \dots$, then $f(z)$ is identically zero.
\end{theorem}
By \eqref{21.02}, $p^{\underline{n}}=0$, for $n+1\geq p$ where $p$ is an arbitrary integer. This implies that $F_0(p)$ converges for any integer $p$, where $F_0(z)$ is given by \eqref{22.02} with $a=0$.
It follows from \eqref{1.7} that
\begin{align}
F_0(p)=\sum_{n=0}^\infty \frac{\Delta^n f(0)}{n!}p^{\underline{n}}=\sum_{n=0}^p \begin{pmatrix}
p\\
n
\end{pmatrix} \Delta^n f(0)=f(p),\quad p=0, 1, 2, \dots .
\end{align}
This gives that $f(p)-F_0(p)=0$, for $p=0, 1, 2, \dots$. By means of Theorem~\ref{ThmCarlson}, we obtain the following Corollary.

\begin{cor}\label{CorCarlson}
Let $f(z)$ be an entire function of order $\rho$ and of type $\tau$. If $\rho<1$, or $\rho=1$ and $\tau<\pi$,
Then $f(z)$ can be represented by $F_0(z)$, i.e.,
\begin{equation}
f(z)=\sum_{n=0}^\infty \frac{\Delta^n f(0)}{n!}z^{\underline{n}}.\label{22.10}
\end{equation}
\end{cor}
Corollary~\ref{CorCarlson} mentions that an entire function $f(z)$ can be represented by a  binomial series under certain growth conditions of $f(z)$. Conversely, Theorem~\ref{Convergence} states that a binomial series gives an entire function with some condition given by its coefficients.

\section{Preliminaries}\label{Preliminaries}
First, we recall the Euler formula for Gamma function, see e.g.,~\cite[Page 257]{M-Thomson1933},~\cite[Page 237]{WW1927}.
\begin{equation}
\Gamma(z)=\lim_{n\to\infty}\frac{(n-1)!\ n^z}{z(z+1)\cdots(z+n-1)}.\label{2.2}
\end{equation}
Putting $-z$ in place of $z$ in \eqref{2.2}, we obtain
\begin{equation}
z^{\underline{n}}=\frac{(-1)^n (n-1)!n^{-z}}{\Gamma(-z)}(1+o(1))=\frac{(-1)^n n!}{\Gamma(-z)n^{1+z}}(1+o(1)),\label{2.3}
\end{equation}
as $n\to\infty$. Let $\{\phi_n\}$ and $\{\psi_n\}$ be positive sequences. We write $\phi_n\sim \psi_n$ if $\phi_n=O(\psi_n)$ and $\psi_n=O(\phi_n)$ as $n\to\infty$, simultaneously.
Since $\Gamma(z)$ has no zeros in $\mathbb C$, we obtain the following lemma by \eqref{2.3}.
\begin{lem}\label{Euler} \enspace For any given $z\in \mathbb C\setminus \mathbb N$, we have
\begin{equation}
|z^{\underline{n}}|\sim \frac{n!}{n^{1+\Re z}},\quad \text{as $n\to\infty$}.\label{2.1}
\end{equation}
\end{lem}
Next we give some estimates on the Stirling numbers of the falling factorial.

\begin{lem}\label{binomial_coefficient}\enspace Let $n\geq2$ be an integer, and $z^{\underline{n}}$ be a falling factorial defined by \eqref{1.3}. Let $\eta_{j, n}$ and $\widetilde{\eta}_{j, n}$ be the Striling numbers defined in \eqref{21.03} and \eqref{21.05}, respectively.
Then we have
\begin{equation}
|\eta_{j,n}|\leq \left(\frac{(n-1)!}{(j-1)!}\right)^2\frac{1}{(n-j)!},\quad
1\leq j\leq n\label{2.23}
\end{equation}
and
    \begin{equation}
|\widetilde{\eta}_{k,n}|\leq \left(\frac{(n-1)!}{(k-1)!}\right)^2\frac{1}{(n-k)!},\quad
1\leq k\leq n.\label{unperbound.eq}
    \end{equation}
\end{lem}

\noindent{\sf Proof of Lemma~\ref{binomial_coefficient}}\quad
We first give the proof of \eqref{2.23} by induction.
Clearly, \eqref{2.23} holds for $n=2$ and $1\leq j\leq2$.
Assume that \eqref{2.23} holds for $n$ and all $1\leq j\leq n$. By means of \eqref{21.04}, we have
\begin{align*}
|\eta_{j,n+1}|&\leq |\eta_{j-1,n}|+n|\eta_{j,n}|\\
&\leq \left(\frac{(n-1)!}{(j-2)!}\right)^2\frac{1}{(n-j+1)!}+n\left(\frac{(n-1)!}{(j-1)!}\right)^2\frac{1}{(n-j)!}\\
&\leq \left(\frac{n!}{(j-1)!}\right)^2\frac{1}{(n-j+1)!}\left(\frac{(j-1)^2+n(n-j+1)}{n^2}\right).
\end{align*}
Combining this and $(j-1)^2+n(n-j+1)-n^2=(j-1)(j-(n+1))\leq0$, we obtain \eqref{2.23} for the case $n+1$ and all $1\leq j\leq n+1$. We have thus proved~\eqref{2.23}.
\vspace{0.25cm}

\noindent
Next we give the proof of \eqref{unperbound.eq}.
Clearly, \eqref{unperbound.eq} holds for $n=2$ and $1\leq k\leq 2$.
Assume that \eqref{unperbound.eq} holds for $n$ and all $1\leq k\leq n$.
By means of \eqref{21.06}, we have
\begin{align*}
|\widetilde{\eta}_{k,n+1}|&\leq |\widetilde{\eta}_{k-1,n}|+k|\widetilde{\eta}_{k,n}|\leq
|\widetilde{\eta}_{k-1,n}|+n|\widetilde{\eta}_{k,n}|\\
&\leq \left(\frac{(n-1)!}{(k-2)!}\right)^2\frac{1}{(n-k+1)!}+n\left(\frac{(n-1)!}{(k-1)!}\right)^2\frac{1}{(n-k)!}\\
&\leq \left(\frac{n!}{(k-1)!}\right)^2\frac{1}{(n-k+1)!}\left(\frac{(k-1)^2+n(n-k+1)}{n^2}\right).
\end{align*}
Combining this and $(k-1)^2+n(n-k+1)-n^2=(k-1)(k-(n+1))\leq0$, we obtain \eqref{unperbound.eq} for the case $n+1$ and all $1\leq k\leq n+1$. We have thus proved Lemma~\ref{binomial_coefficient}.
\quad $\square$

\section{Proof of Theorem~\ref{Convergence}}\label{SConvergence}

Before we start the proof of Theorem~\ref{Convergence}, we refer the Stirling formula
\begin{equation}
n!\sim n^n\frac{\sqrt{2\pi n}}{e^n},\quad \text{as $n\to\infty$},\label{2.4}
\end{equation}
which is useful to estimate asymptotic behavior of $\{a_n\}$ in \eqref{1.5}, see e.g.,~\cite[Page 58]{Titchmarsh1939},~\cite[Page 253]{WW1927}. 
\vspace{0.25cm}

\noindent{\sf Proof of Theorem~\ref{Convergence}}\quad First we show that the formal series $Y(z)$ converges to an entire function $y(z)$ and $\sigma(y)\leq \chi(\{a_n\})$.
By \eqref{3.1} for any $\varepsilon>0$, we have
\begin{equation*}
\frac{n\log n}{-\log |a_n|}\leq \chi(\{a_n\})+\varepsilon,
\end{equation*}
for $n$ sufficiently large. We choose $\varepsilon$ small enough to satisfy $\chi(\{a_n\})+\varepsilon<1$. For the convenience, we set $\gamma=1/(\chi(\{a_n\})+\varepsilon)>1$. Then it holds $|a_n|<n^{-\gamma n}$ for $n$ sufficiently large.
By definition, $z^{\underline{n}}$ vanishes for large $n$ if $z$ is a positive integer, which implies that $Y(z)$ converges for such $z$.
In what follows, we assume that $z$ is not a positive integer.
By means of Lemma~\ref{Euler} and \eqref{2.4}, we have
\begin{equation}
|a_n||z^{\underline{n}}|\leq n^{-\gamma n}\frac{n!}{n^{1+\Re z}}\sim \frac{\sqrt{2\pi n}}{e^n}\frac{1}{n^{1+\Re z}}n^{(1-\gamma) n}.\label{3.2}
\end{equation}
Applying d'Alembert's ratio test, see e.g.,~\cite[Page 22]{WW1927}, to the right hand side of \eqref{3.2}, we see that $Y(z)$ converges to an entire function $y(z)$ uniformly on every compact subset in $\mathbb C$.
\vspace{0.25cm}

Below we show $\rho(y)$ coincides with $\chi(\{a_n\})$. We write
\begin{equation}
y(z)=\sum_{n=0}^\infty b_n z^n,\quad b_n\in\mathbb C,\quad n=0,1, 2, \dots.\label{3.3}
\end{equation}
By \eqref{1.5} and \eqref{21.03}, we have
\begin{equation}
y(z)=\sum_{n=0}^\infty a_n z^{\underline{n}}=\sum_{n=0}^\infty a_n \left(\sum_{j=0}^n \eta_{j, n}z^j \right)=\sum_{n=0}^\infty \left(\sum_{k=n}^\infty a_k \eta_{n, k} \right)z^n,\label{3.4}
\end{equation}
where $\eta_{n, n}=1$, $\eta_{n,0}=0$. It follows from \eqref{3.3} and \eqref{3.4}, $b_n=\sum_{k=n}^\infty a_k \eta_{n, k} $.
By means of~\eqref{2.23} in Lemma~\ref{binomial_coefficient} and \eqref{2.4}, we have
\begin{align*}
|b_n|&\leq \sum_{k=n}^\infty |a_k|| \eta_{n, k}| \leq K_1\sum_{k=n}^\infty k^{-\gamma k}\left(\frac{(k-1)!}{(n-1)!}\right)^2\frac{1}{(k-n)!}\\
&\leq K_2\frac{n^2}{(n!)^2}\sum_{k=n}^\infty k^{-\gamma k}\frac{(k!)^2}{k^2}\frac{1}{(k-n)!}\\
&\leq K_3\frac{n^2}{(n!)^2}\sum_{k=n}^\infty k^{-\gamma k}\frac{1}{k^2}\left(\frac{k^{2k}k}{e^{2k}}\right)\left(k^{-k}e^k\frac{k^n}{\sqrt{k}}\right)\\
&\leq K_3\frac{n^2}{(n!)^2}\sum_{k=n}^\infty k^{(1-\gamma) k}\frac{k^n}{k\sqrt{k}e^k},
\end{align*}
where $K_1$, $K_2$ and $K_3$ are some positive constants. Since $\gamma>1$ and $k^n/e^k$ decreases for $k\geq n$, we have
\begin{align*}
|b_n|&\leq K_4\frac{n^2}{(n!)^2}\sum_{k=n}^\infty n^{(1-\gamma) k}\frac{n^n}{e^n} \leq K_5\frac{n^2}{(n!)^2} n^{(1-\gamma) n}\frac{n^n}{e^n},
\end{align*}
where $K_4$ and $K_5$ are some positive constants. Using \eqref{2.4} again, we obtain
\begin{equation}
|b_n|\leq K_6 n^2\left(\frac{e^{n}}{n^n \sqrt{n}}\right)^2n^{(1-\gamma) n}\frac{n^n}{e^n}\leq K_6 \frac{ne^n}{n^{\gamma n}}.\label{3.5}
\end{equation}
By the Lindel\"of--Pringsheim theorem and \eqref{3.5}, we obtain $\rho(y)\leq \chi(\{a_n\})$.

Secondly we show $\chi(\{a_n\})\leq \rho(y)$.
Using the Lindel\"of--Pringsheim theorem, we obtain $|b_n|<n^{-\tilde\gamma n}$ for $n$ sufficiently large, where
$\tilde\gamma=1/(\rho(y)+\varepsilon)>1$, if we choose $\varepsilon$ small enough.
Therefore, it follows
$$
y(z)=\sum_{n=0}^\infty b_n z^{n}=\sum_{n=0}^\infty b_n \left(\sum_{k=0}^n \widetilde{\eta}_{k, n}z^{\underline{k}} \right)=\sum_{n=0}^\infty \left(\sum_{k=n}^\infty b_k \widetilde{\eta}_{n, k} \right)z^{\underline{n}},
$$
where $\widetilde{\eta}_{n, n}=1$, $\widetilde{\eta}_{n,0}=0$. It follows that $a_n=\sum_{k=n}^\infty b_k \widetilde{\eta}_{n, k} $.
By the same discussion above with~\eqref{unperbound.eq} in Lemma~\ref{binomial_coefficient}, we have
$$
|a_n|\leq \sum_{k=n}^\infty |b_k|| \widetilde{\eta}_{n, k}| \leq K \frac{ne^n}{n^{\tilde\gamma n}},
$$
for some constant $K>0$. By \eqref{3.1}, we obtain $\chi(\{a_n\})\leq \rho(y)$.
Thus, Theorem~\ref{Convergence} is proved.\quad $\square$\vspace{0.25cm}

If $a_n=K a^n/n!$ in \eqref{1.5}, where $K$ is a positive constant and $a$ is a nonzero constant, then we have $\chi(\{a_n\})=1$ from \eqref{3.1}.
From Lemma~\ref{Euler}, we have $|a_n||z^{\underline{n}}|\leq \tilde K|a^n|/n^{1+\Re z}$ for some constant $\tilde K$.
Using the similar arguments in the proof of Theorem~\ref{Convergence}, we obtain the following corollary.

\begin{cor} \label{Convergencechi=1}
Suppose that $|a_n|\leq K/n!$ for some positive constant $K$. Then the formal series $Y(z)$ given by \eqref{1.5} converges to $y(z)$ uniformly on every compact subset in the right half plane $\Re z>0$.
\end{cor}
As an application of Theorem~\ref{Convergence} and Corollary~\ref{Convergencechi=1} to difference equation, we have the following corollary.

\begin{cor} \label{ConvergenceCor}
Let $k\geq1$ be an integer and $R(z,y_1, \dots, y_k)$ be a rational function in $z$ and $y_1$, \dots, $y_k$.
Suppose that the difference equation
\begin{equation}\label{3.6}
y(z)=R(z,y(z+1), \dots, y(z+k))
\end{equation}
has a formal solution $Y(z)$ of the form \eqref{1.5} with $\chi(\{a_n\})<1$ or $|a_n|\leq K/n!$ for some positive constant $K$. Then \eqref{3.6} possesses a meromorphic solution.
\end{cor}

\noindent{\sf Proof of Corollary~\ref{ConvergenceCor}}\quad If $\chi(\{a_n\})<1$, then by Theorem~\ref{Convergence} the formal series $Y(z)$ given by \eqref{1.5} converges to $y(z)$ uniformly on every compact subset in $\mathbb C$.
We consider the case $|a_n|\leq K/n!$ for some positive constant $K$. By means of Corollary~\ref{Convergencechi=1}, $Y(z)$ converges at least in the right half complex plane.
This solution can be analytically continued to the whole complex plane by using \eqref{3.6}. \quad $\square$

\section{Examples}\label{examples}
Let $k$ be a positive integer.
We consider a difference equation
\begin{equation}
\Xi(z,y,\Delta y,\dots, \Delta^{k}y)=0,\label{1.1}
\end{equation}
where $\Xi(x,x_0,x_1,\ldots,x_{k})$ is a polynomial in $x,x_0,x_1,\ldots,x_{k}$.
By means of \eqref{1.6}, it is possible that we write \eqref{1.1} using shifts of $y$ as
\begin{equation}
\tilde\Xi(z,y(z),y(z+1),\dots, y(z+k))=0, \label{1.8}
\end{equation}
where $\Xi(x,x_0,x_1,\ldots,x_{k})$ is a polynomial in $x,x_0,x_1,\ldots,x_{k}$.
Conversely, by \eqref{1.7}, the equation \eqref{1.8} can be written in the form \eqref{1.1}.\vspace{0.25cm}

Before we give examples of binomial series defined by some difference equations, we mention a elemental lemma below, which is useful to manage constructional recurrence relations.

\begin{lem}\label{lem z-times} Let $Y(z)$ be a formal series given by \eqref{1.5}.
Then
\begin{align*}
&zY(z)=\sum_{n=1}^\infty (na_n+a_{n-1})z^{\underline{n}},\quad z\Delta Y(z)=\sum_{n=1}^\infty n\left((n+1)a_{n+1}+a_{n}\right)z^{\underline{n}}\\
&\text{ and }\quad z\Delta^2 Y(z)=\sum_{n=1}^\infty n(n+1)\left((n+2)a_{n+2}+a_{n+1}\right)z^{\underline{n}}.\nn
\end{align*}
\end{lem}

\noindent{\sf Proof of Lemma~\ref{lem z-times}}\quad
Using \eqref{1.5} and \eqref{21.02}, we have
\begin{align*}
(z+1)Y(z)&=\sum_{n=0}^\infty a_n(z+1)^{\underline{n+1}}=\sum_{n=1}^\infty a_{n-1}(z+1)^{\underline{n}}=\sum_{n=1}^\infty a_{n-1}(\Delta z^{\underline{n}}+z^{\underline{n}})\\
&=\sum_{n=1}^\infty a_{n-1}(n z^{\underline{n-1}}+z^{\underline{n}}).
\end{align*}
This gives
\begin{align*}
zY(z)&=(z+1)Y(z)-Y(z)=\sum_{n=1}^\infty a_{n-1}(n z^{\underline{n-1}}+z^{\underline{n}})-\sum_{n=0}^\infty a_nz^{\underline{n}}\\
&=\sum_{n=1}^\infty n a_{n-1}z^{\underline{n-1}}+\sum_{n=1}^\infty
(a_{n-1}-a_n)z^{\underline{n}}-a_0\\
&=\sum_{n=1}^\infty (n+1) a_{n}z^{\underline{n}}+\sum_{n=1}^\infty
(a_{n-1}-a_n)z^{\underline{n}},
\end{align*}
the assertion of the first equality in Lemma~\ref{lem z-times} follows. Similarly, we obtain the second and the third equalities in Lemma~\ref{lem z-times}, respectively\quad $\square$

\begin{rem}\label{general z-times}\enspace {\rm
For any positive integer $k$, we have
$$
z\Delta^k Y(z)=\sum_{n=1}^\infty n(n+1)\cdots (n+k-1)\left((n+k)a_{n+k}+a_{n+k-1}\right)z^{\underline{n}}.
$$
We skip giving a proof of equality just above here. The first equality in Lemma~\ref{lem z-times} implies that $z^k Y(z)=z(z^{k-1}Y(z))$ can be represented by a binomial series
for any integer $k$ recursively. For example,
\begin{equation}
z^2Y(z)=(a_0+a_1)z^{\underline{1}}+\sum_{n=2}^\infty\left(n^2a_{n}+(2n-1)a_{n-1}+a_{n-2}\right)z^{\underline{n}},\label{2.31}
\end{equation}
which is confirmed by elementally computations with $zz^{\underline{n}}=z^{\underline{n+1}}+nz^{\underline{n}}$ and $z^2z^{\underline{n}}=z^{\underline{n+2}}+(2n+1)z^{\underline{n+1}}+n^2z^{\underline{n}}$.
}
\end{rem}

\begin{rem}\label{Periodic}\enspace {\rm Let $\pi$ be a periodic function with period $1$ including a constant. Clearly, for any function $f$, we have $\Delta (\pi(z)f(z))=\pi(z)\Delta f(z)$. We proceed to show that $z^{\underline{n}}$ is linearly independent over periodic function field, that is, if $n_1, \ldots, n_m$, $m\geq2$ are distinct and satisfy
    \begin{equation}
    \pi_1 z^{\underline{n_1}}+\pi_2 z^{\underline{n_2}}+\cdots+\pi_m z^{\underline{n_m}}=0,\label{1.51}
    \end{equation}
where $\pi_j~(1\leq j\leq m)$ are periodic function of period 1, then $\pi_1=\pi_2=\cdots=\pi_m=0$. The proof heavily depends on \eqref{21.02}. Let us assume that $n_1<n_2<\cdots<n_m$, it yields that by \eqref{1.51}
    $$
    n_1\pi_1 z^{\underline{n_1-1}}+n_2\pi_2 z^{\underline{n_2-1}}+\cdots+n_m\pi_m z^{\underline{n_m-1}}=0.
    $$
We continue this way $n_m$ times. Then we have
    $
    (n_m)!\pi_m=0,
    $
which implies $\pi_m=0$. Similarly, we have $\pi_1=\cdots=\pi_m=0$.
}
\end{rem}

\begin{ex}\label{firstorder}{\rm
Let $\lambda(\ne0,\ -1)$, $|\lambda|\leq1$, be a constant. Consider a difference equation
\begin{equation}
\Delta y(z)=\lambda y(z).\label{9.1}
\end{equation}
We see that a formal solution of \eqref{9.1} is decided by $(n+1)a_{n+1}=\lambda a_n$.
This gives that $a_n=a_0\lambda^n/n!$ with an arbitrary constant $a_0$.
Using \eqref{3.1}, we see that $\chi(\{a_0\lambda^n/n!\})=1$.
By means of Corollary~\ref{ConvergenceCor} and the assumption of $\lambda$, we know that \eqref{9.1} has a meromorphic solution.
In fact, $Y(z)=a_0\sum_{n=0}^\infty (\lambda^n/n!)z^{\underline{n}}$ converges to $a_0(1+\lambda)^z$ at least in the right half complex plane.
This solution can be analytically continued to the whole complex plane by using \eqref{9.1}.
}
\end{ex}

\begin{ex}\label{firstorder2}{\rm
Consider a difference equation
\begin{equation}
\Delta y(z)=(z-1)y(z).\label{8.1}
\end{equation}
We consider a formal solution of \eqref{8.1} of the form \eqref{1.5}.
By \eqref{8.1} and Lemma~\ref{lem z-times}, we have
\begin{align*}
\sum_{n=1}^\infty (n+1)a_{n+1}z^{\underline{n}}=\sum_{n=1}^\infty (na_{n}+a_{n-1})z^{\underline{n}}-\sum_{n=0}^\infty a_{n}z^{\underline{n}},
\end{align*}
which gives $a_1+a_0=0$ and $A_{n+1}=A_{n}$, where $A_n=na_{n}+a_{n-1}$, $n\geq1$.
This gives that $a_n=a_0(-1)^n/n!$ with an arbitrary constant $a_0$.
Using \eqref{3.1}, we compute $\chi(\{a_0(-1)^n/n!\})=1$.
By means of Corollary~\ref{ConvergenceCor}, we see that \eqref{8.1} has a meromorphic solution.
In fact,  $Y(z)=a_0\sum_{n=0}^\infty \left((-1)^n/n!\right)z^{\underline{n}}$ converges to $a_0\Gamma(z)$ in the right half complex plane.
This solution can be meromorphically continued to the whole complex plane by using \eqref{8.1}. We also see that this solution admits simple poles at $0, -1, -2, -3,\dots$ by \eqref{8.1}.
}
\end{ex}

\begin{ex}\label{secondorder2}{\rm
Consider a difference equation of second order
\begin{equation}
(4z+6)\Delta^2 y(z)+3\Delta y(z)+y(z)=0.\label{9.11}
\end{equation}
We consider a formal solution to \eqref{9.11}. Using \eqref{21.02} and Lemma~\ref{lem z-times}, we have
\begin{align*}
4&\sum_{n=1}^\infty n(n+1)((n+2)a_{n+2}+a_{n+1})z^{\underline{n}}+
6\sum_{n=0}^\infty (n+1)(n+2)a_{n+2}z^{\underline{n}}\\
&+3\sum_{n=0}^\infty (n+1)a_{n+1}z^{\underline{n}}+\sum_{n=0}^\infty a_{n}z^{\underline{n}}=0,
\end{align*}
which implies
\begin{equation*}
2(n+2)(n+1)(2n+3)a_{n+2}+(n+1)(4n+3)a_{n+1}+a_n=0,\quad n\geq0.
\end{equation*}
Setting $A_n=2(n+1)(2n+1)a_{n+1}+a_{n}$, $n\geq0$, we have
\begin{equation}
(n+1)A_{n+1}+A_n=0,\quad n\geq0.\label{9.15}
\end{equation}
Further, setting $2a_1+a_0=0$, we obtain $A_{0}=0$, and hence by \eqref{9.15}, $A_{n}=0$, $n\geq0$, namely
\begin{equation}
2a_1+a_0=0,\quad 2(n+1)(2n+1)a_{n+1}+a_{n}=0,\quad n\geq0.\label{9.16}
\end{equation}
We consider the convergence of the formal series \eqref{1.5} with the coefficients
$a_n=(-1)^n/(2n)!$, $n\geq 0$. Compute the quantity $\chi(\{a_n\})$ using \eqref{3.1} and \eqref{2.4},
\begin{equation*}
\chi(\{a_n\})=\limsup_{n\to\infty}\frac{n\log n}{\log|(2n)!|}=\frac{1}{2},
\end{equation*}
which concludes that the formal series \eqref{1.5} determined by \eqref{9.16} converges  uniformly on every compact subset in $\mathbb C$ by Theorem~\ref{Convergence}.
}
\end{ex}
\section{Growth of solutions of difference equations}\label{Growth}

In this section, we discuss the order of growth of solutions of difference equations.
First we recall a result of linear difference equations of the form
\begin{equation}
 a_p(z) \Delta^p y(z) + \cdots + a_1(z) \Delta y(z) + a_0(z) y(z) = 0,
\label{4.1}
\end{equation}
where $a_j(z)$ are polynomials of degree $A_j$, $j=0, 1, \dots, p$ and $a_p(z)\not\equiv0$.
Write ${\mathfrak N}_j = \{(x,y) \; ; \; x \geq j, \; y \leq A_{p-j} - (p-j) \}$ for $0 \leq j \leq p.$ The {\it Newton polygon} for \eqref{4.1} is defined as the convex hull of ${\mathfrak N} = \bigcup_{j=0}^p {\mathfrak N}_j.$ The order of growth of entire solutions has been investigated in connection with Newton polygon.
Theorem~\ref{CF} below was showed under the condition $\rho(y)<1/2$ in~\cite[Theorem~1.1]{IY2004}. Chiang and Feng improved the condition on the order of growth using the different method.
\begin{theorem}\label{CF}{\rm ~\cite[Theorem~~4]{CF2016}}\enspace
Let $y$ be an entire solution of order of growth $\rho(y)<1$.
Then the order of growth $\rho(y)$ is a rational number which can be determined from a gradient of the corresponding Newton polygon of \eqref{4.1}.
In particular,
\begin{equation}
\log M(r,y) = Lr^{\rho(y)}(1+o(1)),
\end{equation}
where $L>0$ and $M(r,y)=\max_{|z|=r}|y(z)|$.
\end{theorem}

\begin{rem}\label{NPE}{\rm
We observed the linear difference equation \eqref{9.11} in Example~\ref{secondorder2}. We showed that formal solution $Y(z)$ of
\eqref{9.11} of the form \eqref{1.5} converges to an entire solution $y(z)$ in the complex plane, and we see the growth order $\rho(y)$ is $1/2$ by Theorem~\ref{Convergence}.
Illustrating the Newton polygon of \eqref{9.11}, we ascertain that only $1/2$ is the possible gradient,
 which gives the order of growth $\rho(y)=1/2$ by Theorem~\ref{CF} if $\rho(y)<1$. }
\end{rem}

In what follows, we consider the order of growth of solutions to the difference Riccati equation
\begin{equation}
f(z+1)=\frac{f(z)+A(z)}{1-f(z)},\label{5.1}
 \end{equation}
where $A(z)$ is a rational function.
It is showed \cite[Theorem~3.1]{I2011} that there is no transcendental meromorphic solution to \eqref{5.1} of order less than 1/2, if \eqref{5.1} possesses a rational solution.
Chen and Shon proved \cite[Theorem~1.1]{CS2015} that if \eqref{5.1} possesses a rational solution, then every transcendental meromorphic solution to \eqref{5.1} is of order of growth at least 1. In \cite[Example~3.1]{I2011}, an example of difference Riccati equation of the form \eqref{5.1} that possesses a rational solution and transcendental meromorphic solutions of order at least~1 was constructed.
However, a question which asks does there exist a transcendental meromorphic solution $f$ to \eqref{5.1} of order $0\leq \rho(f)<1$ is generally open. We construct an example of transcendental meromorphic function of order $1/2$ satisfying a difference Riccati equation. We need a lemma below.

\begin{lem}\label{NPlinear} Let $a(\ne0)$, $b$, $c$ be complex constants.
Suppose that a linear difference equation
\begin{equation}
(az+b)\Delta^2y(z)+c\Delta y(z)+y(z)=0\label{5.2}
\end{equation}
possesses a transcendental entire solution $y$ of order $\rho(y)<1$.
Then
\begin{equation}
f(z)=\frac{2(az-a+b)}{2az-2a+2b-c}\left(-\frac{\Delta y(z)}{y(z)}\right)+\frac{c}{2az-2a+2b-c}\label{5.3}
\end{equation}
satisfies a difference Riccati equation of the form \eqref{5.1} with
\begin{equation}
A(z)=\frac{4az-4a+4b+2ac-c^2}{(2az+2b-c)(2az+2b-2a-c)}\label{5.4}
\end{equation}
and $\rho(f)=1/2$.
\end{lem}

\noindent{\sf Proof of Lemma~\ref{NPlinear} }\quad We adopt the method of transformation in \cite[Page 111]{I2017}.
Set $g(z)=-\Delta y(z)/y(z)$.
Then we have $\Delta^2y=-y(z)\Delta g(z)+g(z+1)g(z)y(z)$. Using these formulas and \eqref{5.2}, we see that $g(z)$ satisfies
\begin{equation}
g(z+1)=\frac{1+(az+b-c)g(z)}{az+b-(az+b)g(z)}.\label{5.5}
 \end{equation}
We define $f(z)$ by \eqref{5.3}. It follows from \eqref{5.3} and \eqref{5.5}, we obtain \eqref{5.1} with \eqref{5.4}.
By the Valiron--Mohon'ko theorem, we have $T(r,f)=T(r,g)+O(\log r)$, as $r\to\infty$, see e.g., \cite{M1971},~\cite[Theorem~2.2.5]{L1993}. We have $T(r,g)\leq T(r,y)+T(r,\overline{y})+O(1)$, as $r\to\infty$, where $\overline{y}(z)=y(z+1)$.
Since we assume that $\rho(y)<1$, we have $T(r,\overline{y})=T(r,y)(1+o(1))$,  as $r\to\infty$, $r\not\in E$, where $E\subset \mathbb R$ is a set of finite logarithmic measure, see, e.g., \cite[Theorem~2.1]{CF2008},~\cite[Theorem~8.1]{HKT2014}. Further, by Theorem~\ref{CF} and the arguments in~\cite[Pages 259--261]{HN1984},~\cite[Remark~1.2]{IY2004}, $T(r,\overline{y})=T(r,y)(1+o(1))$ holds for all $r$ sufficiently large, and hence, we obtain $\rho(f)=\rho(g)\leq \rho(y)=1/2$.
By \cite[Theorem~4.2]{GO2008}, we have $\rho(y)=\lambda(y)$, where $\lambda(y)$ is the exponential convergence for zero sequence of $y$. We assert that the number of zeros $z$ of $y$ such that $y(z)=y(z+1)=0$ is finite. If we assume the contrary, there exists $z_0$ satisfying $y(z_0+j)=0$ for all $j\geq1$ by \eqref{5.2}, which implies that $\rho(y)\geq\lambda(y)\geq1$, a contradiction. This gives that $1/2=\lambda(y)=\lambda(1/f)\leq \rho(f)$. Hence we obtain that $\rho(f)=1/2$.\quad $\square$

\begin{prop}\label{Riccati}
There exists a transcendental meromorphic function $f$ of order $\rho(f)=1/2$ satisfying a difference Riccati equation of the form \eqref{5.1}.
\end{prop}\quad

\noindent{\sf Proof of Proposition~\ref{Riccati}}\quad As observed in Example~\ref{secondorder2} and Remark~\ref{NPE} that the second order difference equation \eqref{9.11} possesses an entire solution of order $1/2$. By means of Lemma~\ref{NPlinear}, a difference Riccati equation \eqref{5.1} with
$$
A(z)=\frac{16z+23}{64z^2+80z+9}
$$
has a transcendental meromorphic solution of order $1/2$. \quad $\square$

\end{document}